\documentclass[a4paper,11pt]{article}
\usepackage{graphicx,bm, color}
\usepackage{amsmath,amsfonts,amssymb}
\usepackage[english]{babel}
\newcommand{\sv}{\, ,}
\newcommand{\p}{\, .}
\newcommand{\dm}{\displaystyle}
\newcommand{\perogni}{\forall \,}
\newcommand{\eq}[1]{(\ref{#1})}
\newcommand{\fr}{\\[5pt]}
\newcommand{\lre}{\xi}
\author{Armando Majorana and Rita Tracin\`{a}
	\fr Department of Mathematics and Computer Science, \\
	University of Catania, Italy}
\title{}
\title{Exact and numerical solutions to a Mindlin microcontinuum model}
\begin{document}
\maketitle
%
%
\abstract{
In this paper we consider a one-dimensional Mindlin model describing linear elastic behaviour of isotropic materials with micro-structural effects.
After introducing the kinetic and the potential energy, we derive a system of equations of 
motion by means of the Euler-Lagrange equations. 
A class of exact solutions is obtained. They have a wave behaviour due to a good property of the potential energy.
Numerical solutions are obtained by using a weighted essentially non-oscillatory finite 
difference scheme coupled by a total variation diminishing Runge-Kutta method.
A comparison between exact and numerical solutions shows the robustness and the accuracy 
of the numerical scheme.
A numerical example of solutions for an inhomogeneous material is also shown.
}
%
%
\\[7pt]
MSC-class: 74B - 74J05 (Primary) 74S20 (Secondary)
%
\baselineskip=16pt
\section{Introduction}
In this paper we study a one-dimensional Mindlin model for longitudinal waves in 
microstructured materials.
It is well known that matter is not continuous but has an internal structure. 
This underlying microstructure could influence profoundly the dynamic thermo-mechanical 
response of materials, as, for instance, when the scale of a deformation is of the order 
of the material microstructural heterogeneities or the characteristic length of the continuum is comparable to the short wavelength of a signal propagating through the material.
An interesting example of the impact of microstructural peculiarities on the dynamic 
response of materials is the wave scattering due to the mismatch in mechanical impedance at the interface between different material phases in layered composites \cite{Zhuang2003}. 
In \cite{Millett2007} it has been experimentally observed that heterogeneities, such as, for instance, microvoids leading to localization of deformation and loss of shear strength, have effect on the dynamic strength of polycrystalline metals.\\
There are essentially two types of models to study microstructured materials: discrete and 
continuum models.
In the microcontinuum theory, the macrostructure and microstructure of the continuum are 
usually separated. 
This leads to the formulation of separate balance laws for each structure.
Mindlin \cite{Mindlin1964} proposed a model, where the two structures are separated, and 
introduced two balance laws: one for the macrostructure and another for the microstructure.
Moreover the Mindlin continuum model also incorporates the inelastic behaviour of the 
material.
In the interesting Ref.~\cite{Berezovski2011}, \cite{Berezovski2016b} and \cite{Engelbrecht2015} general one-dimensional Mindlin-type microstructure models are discussed.\\
In this work we transform the equations of a one-dimensional Mindlin model in a particular set of hyperbolic partial differential equations, where it is clear how to impose correct boundary conditions. This point is important in the numerical simulations. 
Moreover we propose a numerical scheme, based on a weighted essentially non-oscillatory finite difference scheme coupled by a total variation diminishing Runge-Kutta method, which furnishes accurate results also for not smooth solutions. \\	
This paper is organized as follows. In Section 2 we derive the one-dimensional wave equations for microstructured materials and we find important inequalities involving the physical parameters, by requiring the potential energy to be a strictly positive definite function.
Section 3 is devoted to look for simple exact solutions and in Section 4 we transform the 
equations in a suitable system of first order partial differential equations, which is 
solved numerically. 
In Section 5 we consider the case when the physical parameters depend on the spatial 
coordinate. This happens when the continuum consists of different materials. Also in 
this case numerical solutions are obtained.
Finally the conclusions are drawn in the last section.
\section{Basic equations}
Following an approach similar to that employed by the authors in \cite{Mindlin1964}, and 
\cite{Engelbrecht2005}, the one-dimensional wave equations for 
microstructured materials are derived starting from the Lagrangian
\begin{equation}
{\cal{L}}={\cal{K}}-{\cal{W}},
\end{equation}
where $\cal{K}$ and $\cal{W}$ are the kinetic and potential energy, respectively. 
We assume that the kinetic energy $\cal{K}$ is given by
\begin{equation}
{\cal{K}} = \frac{1}{2} \rho u_t^2 + \frac{1}{2} I_{\mu} \chi_t^2 \sv
\label{kinetic}
\end{equation}
where $u$ is macroscopic displacement, and $\chi$ the microdeformation.
The constant positive parameters $\rho$ and $I_{\mu}$ are  the macroscopic density and the micro-inertia, respectively. As usually the subscript $t$ or $x$ denotes the first partial 
derivative with respect to the time $t$ or the space coordinate $x$. 
According to the theory of the elasticity, the potential energy $\cal{W}$ is a function 
of the variables $u_x$, $\chi$, $\chi_x$. \\
The corresponding Euler-Lagrange equations have the general form
\begin{align}
&
\left( \dfrac{\partial \cal{L}}{\partial u_{t}} \right)_{\! t} +
\left( \dfrac{\partial \cal{L}}{\partial u_{x}} \right)_{\! \! x} = 0 \sv
\label{EL1} 
\fr
&
\left( \dfrac{\partial \cal{L}}{\partial \chi_{t}} \right)_{\! t} +
\left( \dfrac{\partial \cal{L}}{\partial \chi_{x}} \right)_{\! \! x} -
\dfrac{\partial \cal{L}}{\partial \chi} = 0 \p
\label{EL2}
\end{align}
We choose the potential energy $\cal{W}$ given by
\begin{equation}
{\cal{W}} = \frac{1}{2} \gamma \, u_x^2 +  A \, u_x \, \chi + \frac{1}{2} \, B \, \chi^2 
+ \frac{1}{2} \, C \, \chi_x^2 \sv
\label{potential}
\end{equation} 
where $\gamma$, $A$, $B$, and $C$ are physical constant parameters of the model.
Now the system \eq{EL1}-\eq{EL2} becomes
\begin{align}
&
\rho \, u_{tt} = \gamma \, u_{xx} + A \, \chi_x  \sv
\label{equ_tt}
\fr
&
I_\mu \, \chi_{tt} = C \, \chi_{xx} - A \, u_x - B \, \chi \p
\label{eqx_tt}
\end{align}
The case $A = 0$ is mathematically interesting, because Eq.~\eq{equ_tt} reduces to 
the classical wave equation for the unknown $u$, provided that $\gamma > 0$, and 
Eq.~\eq{eqx_tt} becomes the telegrapher's equation, where $\chi$ is the unknown. We do not 
consider this special case.  \\
It is also possible to derive a single fourth-order partial differential equations from the system \eq{equ_tt}-\eq{eqx_tt} (see also Ref.~\cite{Metrikine2006}). 
In fact Eq.~\eq{equ_tt} gives
$$
\chi_{x} = \dfrac{1}{A} \left[ \rho \, u_{tt} - \gamma \, u_{xx} \right]
$$
Now, if we differentiate Eq.~\eq{eqx_tt} with respect to the variable $x$ and use the 
previous equation, then it is a simple matter to obtain the following equation
$$
I_{\mu}  \, \rho \, u_{tttt} - \left( \gamma \,  I_{\mu} + \rho \, C \right) u_{ttxx} +
\gamma \, C \, u_{xxxx} + \left( A^{2} - B \, \gamma \right) u_{xx} + B \, \rho \, u_{tt}
 = 0 \p
$$
In this paper we do not study this equation, but we consider the system 
\eq{equ_tt}-\eq{eqx_tt}.

The main assumption of this paper is the following
\begin{center}
\emph{ the potential energy} 
${\cal{W}} = {\cal{W}}( u_{x}, \chi, \chi_{x} )$ 
\emph{is a strictly positive definite function.}
\end{center}
This implies some inequalities involving the physical parameters.
We have immediately
\begin{align}
&
{\cal{W}}( u_{x}, 0, 0 ) = \frac{1}{2} \, \gamma \, u_{x}^{2} > 0
\quad \perogni u_{x} \neq 0 \Rightarrow \gamma > 0 \sv
\label{dis1W}
\fr
&
{\cal{W}}( 0, \chi, 0 ) = \frac{1}{2} \, B \, \chi^{2} > 0
\quad \perogni \chi \neq 0 \Rightarrow B > 0 \sv
\label{dis2W}
\fr
&
{\cal{W}}( 0, 0, \chi_{x} ) = \frac{1}{2} \, C \, \chi_{x}^{2} > 0
\quad \perogni \chi_{x} \neq 0 \Rightarrow C  > 0 \p
\label{dis3W}
\end{align}
Since
$ 
\dm {\cal{W}}( u_{x}, \chi, \chi_{x} ) = {\cal{W}}( u_{x}, \chi, 0 ) +
{\cal{W}}( 0, 0, \chi_{x} )
$,
then, thanks to Eq.~\eq{dis3W}, ${\cal{W}}$ is a strictly positive definite function if and only if
$$
{\cal{W}}( u_{x}, \chi, 0 ) =
\frac{1}{2} \, \gamma \, u_x^2 + A \, u_x \, \chi + \frac{1}{2} \, B \, \chi^2 
$$
is a strictly positive definite function.
This holds if and only if the discriminant is negative.
This gives the last condition
\begin{equation}
\gamma \, B  - A^{2} > 0 \p
\label{dis4W}
\end{equation}
\section{Explicit simple solutions}
In this section we are interested to derive some explicit solutions to the one-dimensional 
wave equations for microstructured materials ~\eq{equ_tt}-\eq{eqx_tt}.
To the best of our knowledge, the complete derivation of the following solutions has never 
been reported, although in \cite{Berezovski2011}, \cite{Berezovski2016}, \cite{Dingreville}, 
\cite{Engelbrecht2005}, \cite{Peets2010} the authors carried out a dispersion analysis 
for the same equations.  
\\
To simplify notation, it is useful to define the parameters
\begin{equation}
a_{1} = \frac{\gamma}{\rho} \sv \quad
a_{2} = \frac{A}{\rho} \sv \quad
a_{3} = \frac{C}{I_\mu} \sv \quad
a_{4} = \frac{A}{I_\mu} \sv \quad
a_{5} = \frac{B}{I_\mu} \p
\end{equation}
Hence Eqs.~\eq{equ_tt}-\eq{eqx_tt} write
\begin{align}
& u_{tt} = a_{1} \, u_{xx} + a_{2} \, \chi_x \sv
\label{au_tt}
\\
& \chi_{tt} = a_{3} \, \chi_{xx} - a_{4} \, u_x - a_{5} \, \chi  \p 
\label{ax_tt}
\end{align}
The inequalities \eq{dis1W}-\eq{dis4W} imply that
\begin{equation}
a_{1} > 0 \sv \quad  a_{3} > 0 \sv \quad  a_{5} > 0  \sv
\label{ineq_a}
\end{equation}
and
\begin{equation}
a_{2} \, a_{4} - a_{1} \, a_{5} = 
\dfrac{A}{\rho} \, \dfrac{A}{I_{\mu}} - \dfrac{\gamma}{\rho} \, \dfrac{B}{I_{\mu}} 
= \dfrac{1}{\rho \, I_{\mu}} \left( A^{2} - \gamma \, B \right) < 0 \p
\label{ineq_b}
\end{equation}
\indent
We look for solutions of Eqs.~\eq{au_tt}-\eq{ax_tt} of the kind
\begin{equation}
u(t,x) = U_{1}(t) \, \sin (\omega \, x) \sv \quad 
\chi(t,x) = X_{1}(t) \, \cos (\omega \, x) \sv
\label{U1X1w} 
\end{equation}
where $\omega$ is a non-zero real parameter.
Using \eq{U1X1w}, Eqs.~\eq{au_tt}-\eq{ax_tt} give the system of second order ordinary 
differential equations for the unknowns $U_{1}$ and $X_{1}$
\begin{align*}
U_{1}''(t) \, \sin (\omega \, x) & = 
- \, a_{1} \, \omega^{2} \, U_{1}(t) \, \sin (\omega \, x) 
- a_{2} \, \omega \, X_{1}(t) \,  \sin (\omega \, x) \sv
\fr
X_{1}''(t) \, \cos (\omega \, x) & = 
- \, a_{3} \, \omega^{2} \, X_{1}(t) \, \cos (\omega \, x) 
- a_{4} \, \omega \, U_{1}(t) \, \cos (\omega \, x) - 
a_{5} \, X_{1}(t) \, \cos (\omega \, x) \p
\end{align*}
It is simplified immediately, and writes
\begin{align}
U_{1}''(t) & = \mbox{} - a_{1} \, \omega^{2} \, U_{1}(t) - a_{2} \, \omega \, X_{1}(t) \sv
\label{U1a}
\fr
X_{1}''(t) & = \mbox{} - a_{3} \, \omega^{2} \, X_{1}(t) -
a_{4} \, \omega \, U_{1}(t) - a_{5} \, X_{1}(t) \p
\label{X1a}
\end{align}
Since $\dm a_{2} = A/\rho$ is a non-zero number, then Eq.~\eq{U1a} gives
\begin{equation}
X_{1}(t) = \dfrac{-1}{a_{2} \, \omega} \left[ 
U_{1}''(t) + a_{1} \, \omega^{2} \, U_{1}(t) \right] ,
\label{eqX1}
\end{equation}
and Eq.~\eq{X1a} becomes
\begin{equation}
U_{1}^{(4)}(t) + \left( a_{1} \, \omega^{2} + a_{3} \, \omega^{2} + 
a_{5} \right) U_{1}''(t) + \omega^{2} \left[ 
a_{1} \left( a_{3} \, \omega^{2} + a_{5} \right) - a_{2} \, a_{4} \right] U_{1}(t) = 0 \p
\label{eqU1}
\end{equation}
Eq.~\eq{eqU1} is a linear homogeneous fourth-order ordinary differential equation with 
constant coefficients and then it can be solved easily.
The characteristic equation is 
\begin{equation}
\lambda^{4} + \left( a_{1} \, \omega^{2} + a_{3} \, \omega^{2} + a_{5} \right) \lambda^{2} +
\omega^{2} \left[ a_{1} \, a_{3} \, \omega^{2} + \left(a_{1} \, a_{5} - a_{2} \, a_{4}  
\right) \right] = 0 \sv
\label{eq_car_a}
\end{equation}
where $\lambda$ is an eigenvalue.
The discriminant of the biquadratic equation \eq{eq_car_a} is
\begin{align*}
\Delta & = \left( a_{1} \, \omega^{2} + a_{3} \, \omega^{2} + a_{5} \right)^{2} - 4 \,
\omega^{2} \left[ a_{1} \, a_{3} \, \omega^{2} + \left(a_{1} \, a_{5} - a_{2} \, a_{4}  
\right) \right] 
\\
& = \left[ (a_{3} - a_{1}) \, \omega^{2} + a_{5}\right]^{2} + 4 \, a_{2} \, a_{4} 
\, \omega^{2} \p
\end{align*}
Since $a_{2} \, a_{4}$ is positive, then also $\Delta$ is positive.
Taking into account \eq{ineq_a} and \eq{ineq_b}, it is evident that all the coefficients of Eq.~\eq{eq_car_a} are positive and therefore Eq.~\eq{eq_car_a} has four purely imaginary roots.
If we denote by $\pm i \, \lre$ and $\pm i \, \eta$ the four roots, then the general 
solution of Eq.~\eq{eqU1} writes
\begin{equation}
U_{1}(t) = k_{1} \, \cos(\lre \, t) + k_{2} \, \sin(\lre \, t) +
k_{3} \, \cos(\eta \, t) + k_{4} \, \sin(\eta \, t) \sv
\end{equation}
where $k_{i}$ $(i=1,2,3,4)$ are arbitrary real numbers.
Eq.~\eq{eqX1} gives the solution $X_{1}(t)$ easily
\begin{equation}
X_{1}(t) = \dfrac{\lre^2 - a_{1} \, \omega^2}{a_{2}\, \omega}
\left[ k_{1} \, \cos(\lre \, t) + k_{2} \, \sin(\lre \, t) \right] +
\dfrac{\eta^2 - a_{1} \, \omega^2}{a_{2}\, \omega}
\left[ k_{3} \, \cos(\eta \, t) + k_{4} \, \sin(\eta \, t) \right] .
\end{equation}
Another set of exact solutions can be obtained, looking for solutions of 
Eqs.~\eq{au_tt}-\eq{ax_tt} of the kind 
\begin{equation}
u(t,x) = U_{2}(t) \, \cos (\omega \, x) \sv \quad 
\chi(t,x) = X_{2}(t) \, \sin (\omega \, x) \p
\label{U2X2w} 
\end{equation}
Using \eq{U2X2w}, Eqs.~\eq{au_tt}-\eq{ax_tt} give the system of second order ordinary 
differential equations for the unknowns $U_{2}$ and $X_{2}$
\begin{align*}
U_{2}''(t) \, \cos (\omega \, x) & = 
\mbox{} - a_{1} \, \omega^{2} \, U_{2}(t) \, \cos (\omega \, x) +
a_{2} \, \omega \, X_{2}(t) \, \cos (\omega \, x) \sv
\fr
X_{2}''(t) \, \sin (\omega \, x) & =
\mbox{} - a_{3} \, \omega^{2} \, X_{2}(t) \, \sin (\omega \, x)  
+ a_{4} \, \omega \, U_{2}(t) \, \sin (\omega \, x) - 
a_{5} \, X_{2}(t) \, \sin (\omega \, x) \p
\end{align*}
It is equivalent to the system
\begin{align*}
U_{2}''(t) & = \mbox{} - 
a_{1} \, \omega^{2} \, U_{2}(t) + a_{2} \, \omega \, X_{2}(t) \sv
\fr
X_{2}''(t) & = \mbox{} - 
a_{3} \, \omega^{2} \, X_{2}(t) + a_{4} \, \omega \, U_{2}(t)  - a_{5} \, X_{2}(t) \sv
\end{align*}
that is similar to the system \eq{U1a}-\eq{X1a}.
So we can use the same procedure to derive a fourth order ordinary differential equation.
Since the characteristic equation coincides with Eq.~\eq{eq_car_a}, we do not give 
further details on the solutions of this equation. \\
We remark that Eqs.~\eq{au_tt}-\eq{ax_tt} are linear and homogeneous; so any linear 
(finite or, under suitable conditions, numerable) combination of solutions is also a 
solution.
\section{Numerical solutions and numerical tests}
The numerical treatment of Eqs.~\eq{au_tt}-\eq{ax_tt} requires some transformations of 
variables in order to derive a suitable system of first order partial differential 
equations. 
As the equations are of hyperbolic type, this step is of fundamental importance to 
achieve good numerical solutions, and to take into account boundary conditions, correctly.
\\
Firstly we introduce the new variables $\alpha$ and $\beta$ defined by
\begin{equation}
\alpha = u_{t} - \sqrt{a_{1}} \, u_{x} \sv \quad
\beta = \chi_{t} - \sqrt{a_{3}} \, \chi_{x} \p
\end{equation}
Now Eq.~\eq{au_tt} becomes
$$
\dfrac{\partial \mbox{ }}{\partial t} \left( \alpha + \sqrt{a_{1}} \, u_{x} \right)
- \sqrt{a_{1}} \, \dfrac{\partial \mbox{ }}{\partial x} \left( u_{t} - \alpha \right) 
- a_{2} \, \chi_x = 0
\quad \Leftrightarrow \quad
\alpha_{t} + \sqrt{a_{1}} \, \alpha_{x} - a_{2} \, \chi_x = 0 \sv
$$
and Eq.~\eq{ax_tt} is
$$
\dfrac{\partial \mbox{ }}{\partial t} \left( \beta + \sqrt{a_{3}} \, \chi_{x} \right)
- \sqrt{a_{3}} \, \dfrac{\partial \mbox{ }}{\partial x} \left( \chi_{t} - \beta \right)
+ a_{4} \, u_x + a_{5} \, \chi = 0 
\quad \Leftrightarrow \quad
\beta_{t} + \sqrt{a_{3}} \, \beta_{x} + a_{4} \, u_x + a_{5} \, \chi = 0 \p
$$
Hence, system \eq{au_tt}-\eq{ax_tt} is equivalent to the set of four partial differential 
equations 
\begin{align}
u_{t} = & \sqrt{a_{1}} \, u_{x} + \alpha \sv
\fr
\chi_{t} = & \sqrt{a_{3}} \, \chi_{x} + \beta \sv 
\fr
\alpha_{t} = & \mbox{} - \sqrt{a_{1}} \, \alpha_{x} + a_{2} \, \chi_{x} \sv
\fr
\beta_{t} = & \mbox{} - \sqrt{a_{3}} \, \beta_{x} - a_{4} \, u_{x} - a_{5} \, \chi \p
\end{align}
At this step it is necessary the change of variables
$$
v = \alpha + c_{1} \, \chi \sv \quad
w = \beta + c_{2} \, u
$$
where $c_{1}$ e $c_{2}$ are real constants to be determined.
The new system writes
\begin{align*}
\dfrac{\partial u}{\partial t} & =
\sqrt{a_{1}} \, \dfrac{\partial u}{\partial x} + v - c_{1} \, \chi \sv
\fr
\dfrac{\partial \chi}{\partial t} & =
\sqrt{a_{3}} \, \dfrac{\partial \chi}{\partial x} + w - c_{2} \, u \sv 
\fr
\dfrac{\partial v}{\partial t} & =
c_{1} \, \dfrac{\partial \chi}{\partial t} - \sqrt{a_{1}} \left( 
\dfrac{\partial v}{\partial x} - c_{1} \, \dfrac{\partial \chi}{\partial x}
\right) + a_{2} \, \dfrac{\partial \chi}{\partial x} \sv
\fr
\dfrac{\partial w}{\partial t} & =
c_{2} \, \dfrac{\partial u}{\partial t} - \sqrt{a_{3}} \left( 
\dfrac{\partial w}{\partial x} - c_{2} \, \dfrac{\partial u}{\partial x}
\right) - a_{4} \, \dfrac{\partial u}{\partial x} - a_{5} \, \chi \sv  
\end{align*}
that is
\begin{align*}
\dfrac{\partial u}{\partial t} & =
\sqrt{a_{1}} \, \dfrac{\partial u}{\partial x} + v - c_{1} \, \chi \sv
\fr
\dfrac{\partial \chi}{\partial t} & =
\sqrt{a_{3}} \, \dfrac{\partial \chi}{\partial x} + w - c_{2} \, u \sv 
\fr
\dfrac{\partial v}{\partial t} & =
c_{1} \left[ \sqrt{a_{3}} \, \dfrac{\partial \chi}{\partial x} + w - c_{2} \, u
\right] - \sqrt{a_{1}} \left( 
\dfrac{\partial v}{\partial x} - c_{1} \, \dfrac{\partial \chi}{\partial x}
\right) + a_{2} \, \dfrac{\partial \chi}{\partial x} \sv
\fr
\dfrac{\partial w}{\partial t} & =
c_{2} \left[ \sqrt{a_{1}} \, \dfrac{\partial u}{\partial x} + v - c_{1} \, \chi
\right] - \sqrt{a_{3}} \left( 
\dfrac{\partial w}{\partial x} - c_{2} \, \dfrac{\partial u}{\partial x}
\right) - a_{4} \, \dfrac{\partial u}{\partial x} - a_{5} \, \chi \p
\end{align*}
Now, if we choose
$$
c_{1} = \dfrac{- \, a_{2}}{\sqrt{a_{1}} + \sqrt{a_{3}}}
\quad \mbox{and} \quad
c_{2} = \dfrac{ a_{4}}{\sqrt{a_{1}} + \sqrt{a_{3}}} \sv
$$
the partial derivative of the unknowns $u$ and $\chi$ disappear in the last two 
equations.
The final result is
\begin{align}
\dfrac{\partial u}{\partial t} & =
\sqrt{a_{1}} \, \dfrac{\partial u}{\partial x} - c_{1} \, \chi + v \sv
\label{eqV1}
\fr
\dfrac{\partial \chi}{\partial t} & =
\sqrt{a_{3}} \, \dfrac{\partial \chi}{\partial x} - c_{2} \, u + w \sv 
\fr
\dfrac{\partial v}{\partial t} & =
\mbox{} - \sqrt{a_{1}} \, \dfrac{\partial v}{\partial x} -
c_{1} \, c_{2} \, u + c_{1} \, w \sv
\fr
\dfrac{\partial w}{\partial t} & =
\mbox{} - \sqrt{a_{3}} \, \dfrac{\partial w}{\partial x} -
\left( c_{1} \, c_{2} + a_{5} \right) \chi + c_{2} \, v \p
\label{eqV4}
\end{align}
Eqs.~\eq{eqV1}-\eq{eqV4} is a system of partial differential equations of hyperbolic type, 
and the new unknowns $u, \chi, v, w$ are the Riemann invariants. 
These equations are similar to advection equations with source terms.
This makes simple a proper imposition of the initial and boundary conditions.
For instance, if $x \in [O, L]$, then we must assign the initial conditions for all the variables, and the boundary conditions
\begin{equation}
u(t,L^{+}) \sv \quad \chi(t,L^{+}) \sv \quad
v(t,0^{-}) \sv \quad w(t,0^{-}) \p
\label{bc}
\end{equation}
In order to find numerical solutions of Eqs.~\eq{eqV1}-\eq{eqV4}, we use a Weighted 
Essentially Non-Oscillatory (WENO) finite difference scheme \cite{Liu1994}, \cite{Shu1998}.
The main advantage of WENO schemes is their capability to achieve arbitrarily high-order 
accuracy in regions where the solution are smooth, while maintaining stable, 
non-oscillatory and sharp discontinuity transitions. In particular the schemes are  
suitable for hyperbolic partial differential equations admitting solutions containing both 
strong discontinuities and complex smooth features. 
WENO schemes approximate spatial partial derivatives by means of suitable finite 
differences, so the original system of partial differential equations is replaced by a set of ordinary differential equations, which are solved by means of a Runge-Kutta method. In 
this paper we employ a third-order Total Variation Diminishing (TVD) Runge-Kutta method 
\cite{Jiang1996}.
These methods guarantee that the total variation of the solution does not increase, so that no new extrema are generated.
\\
We consider two simple explicit solutions and make a comparison with the corresponding 
numerical solutions obtained using WENO scheme. The spatial domain of the solutions is the 
interval $(0,1)$.
We assume periodic boundary conditions in the simulations.
The numerical physical parameters, with arbitrary units, used in the simulations are
$$
\rho = 1 \sv \quad I_\mu = 1 \sv \quad
\gamma = 0.99 \sv \quad A = - \, 0.01 \sv \quad B = 10 \sv \quad C = 1 \sv
$$
which are of the same order of the parameters employed in Ref.~\cite{Berezovski2016}.
We point out that the physical parameters are not of same order. This implies that, in general, the solutions are not very smooth.
We use different partitions of the interval $(0,1)$ and we denote by $N$ the number of the 
cells. So the spatial step $\Delta x$ is given by $1/N$, because the length of the interval 
is equal to one. 
We measure the difference between exact and numerical solutions by means of the formula
$$
\mbox{err}(u) = 
\max_{t, x} \dfrac{\left| u(t,x) - \hat{u}(t,x)\right|}{1 + \left|u(t,x) \right| } ,
$$
where $\hat{u}(t,x)$ denotes the numerical solution for the unknown $u$ at 
time $t$ and position $x$; 
of course, only the grid points, in time and space, are used to evaluate the maximum.
An analogous formula is employed for the unknown $\chi$.
\\
{\bf Numerical test A} \\
The parameter $\omega$ of first solution (case A) of type \eq{U1X1w} is equal to $2 \, \pi$, and $k_{i} = 1$ with $(i=1,2,3,4)$. In our simulations, we choose $[0, 10]$ as interval for time integration.
The Table 1 shows the errors between exact and numerical solutions in the case A.
\begin{center}
\begin{table}[ht]
\centering
\caption{Errors for the test problem A.}
\begin{tabular}{|r||c|c|c|c|c|}
\hline
$N$ & $128$ & $256$ & $512$ & $1024$ & $2048$    \\
\hline
err($u$)    & $3.663 \times 10^{-5}$ & $4.471 \times 10^{-6}$ & $5.561 \times 10^{-7}$ 
& $6.947 \times 10^{-8}$ & $8.684 \times 10^{-9}$
\\
\hline	
err($\chi$) & $1.749 \times 10^{-4}$ & $1.259 \times 10^{-5}$ & $1.132 \times 10^{-6}$ 
& $1.184 \times 10^{-7}$ & $1.479 \times 10^{-8}$
\\
\hline
\end{tabular}
\end{table}
\end{center}
\begin{figure}[ht]
\includegraphics[width=0.49\textwidth]{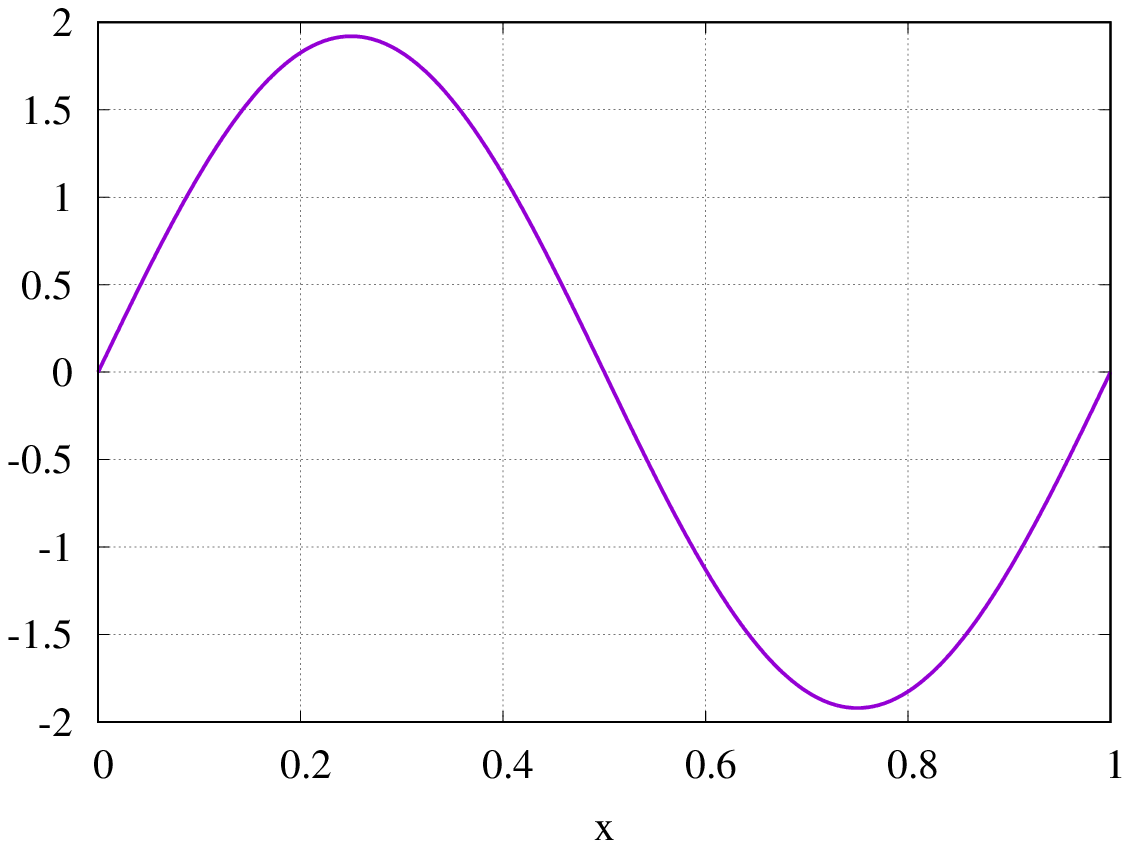} 
\includegraphics[width=0.49\textwidth]{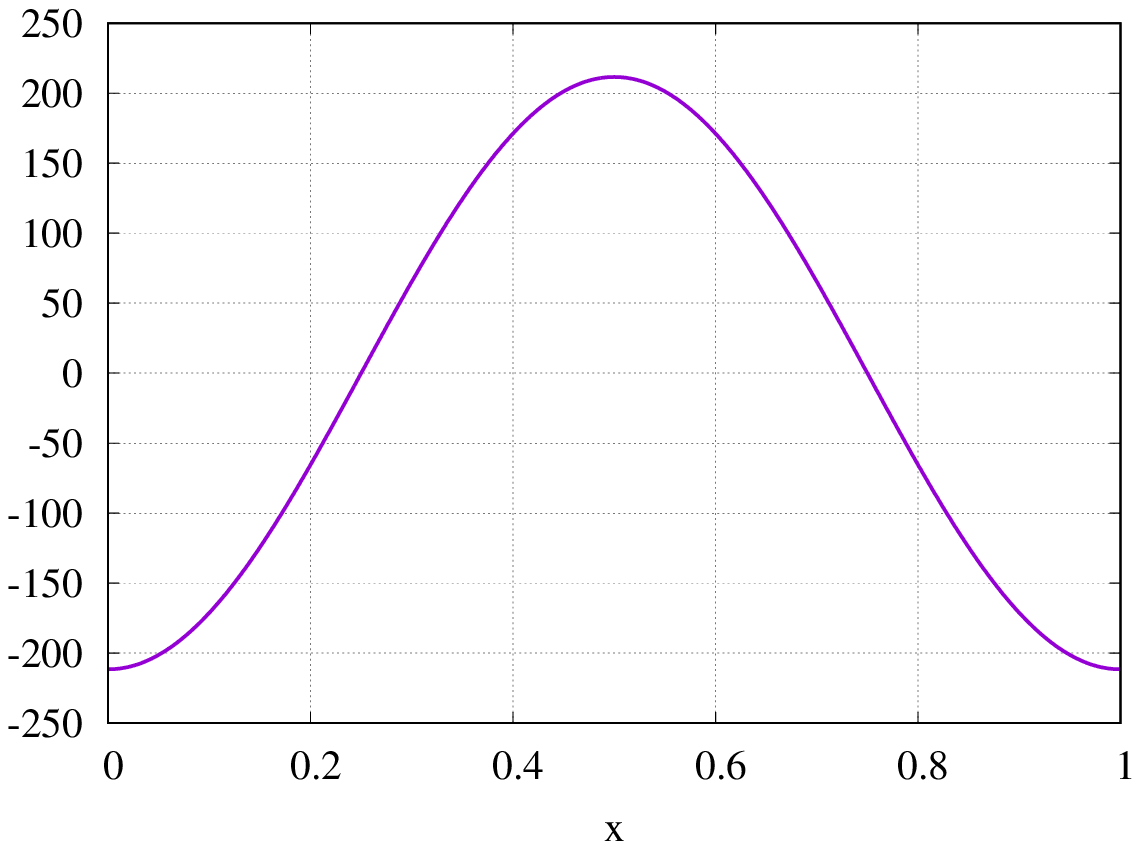} 
\caption{$u(t,x)$ (left) and $\chi(t,x)$ (right) at time $t =10$ in the case A}
\end{figure}
{\bf Numerical test B} \\
For the second test (case B), we use the sum of two solutions of kind \eq{U1X1w} by choosing two values ($2 \, \pi$ and $4 \, \pi$) for the parameter $\omega$. The integration constants $k_{i}$ and the domain parameters are the same as in the case A.
\begin{figure}[ht]
\includegraphics[width=0.49\textwidth]{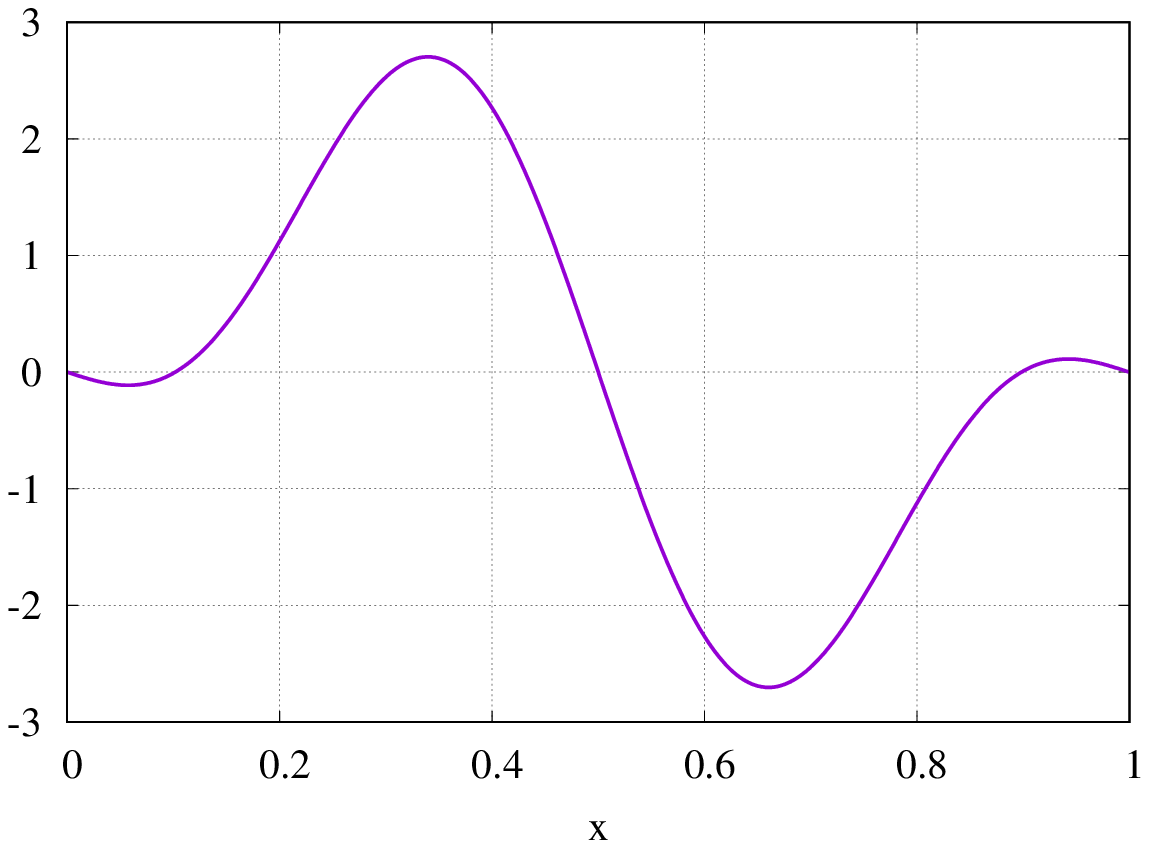} 
\includegraphics[width=0.49\textwidth]{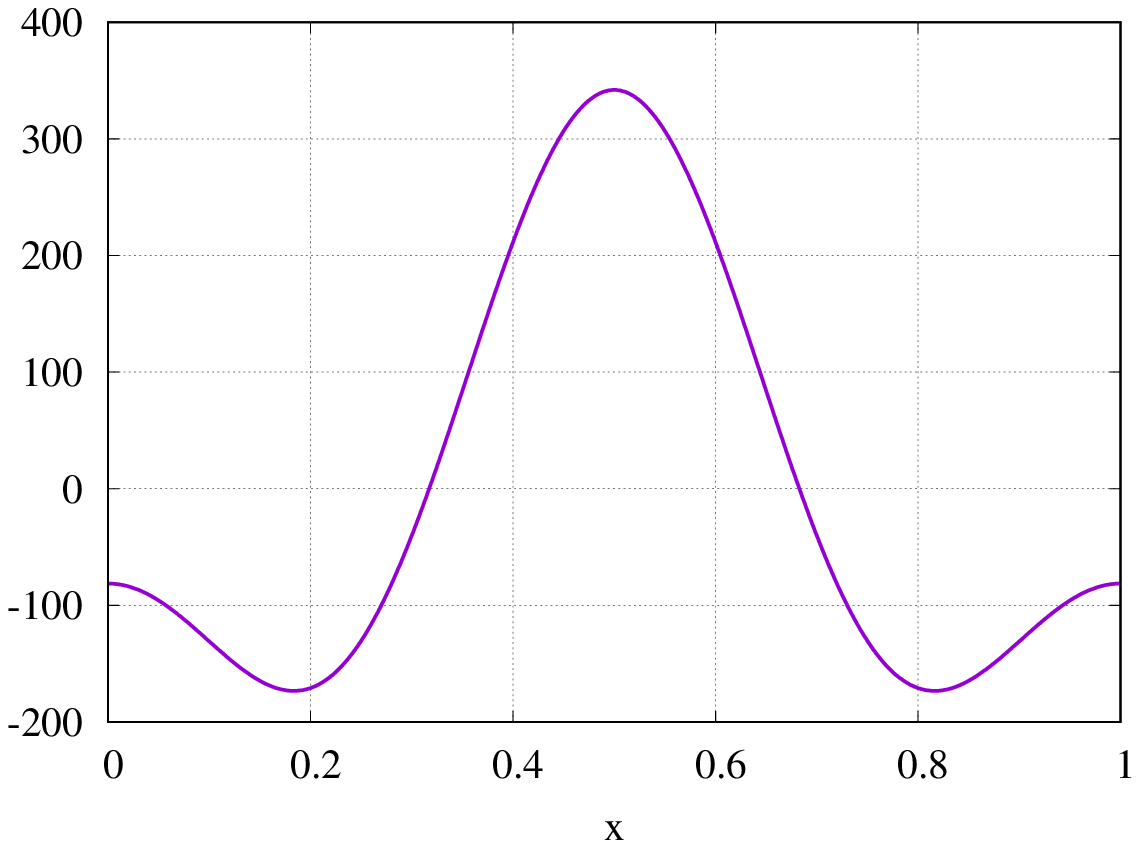} 
\caption{$u(t,x)$ (left) and $\chi(t,x)$ (right) at time $t =10$ in the case B}
\end{figure}
The Table 2 shows the errors between exact and numerical solutions in the case B.
\begin{center}
\begin{table}[!ht]
\centering
\caption{Errors for the test problem B.}
\begin{tabular}{|r||c|c|c|c|c|}
\hline
$N$ & $128$ & $256$ & $512$ & $1024$ & $2048$    \\
\hline
err($u$)    & $1.151 \times 10^{-3}$ & $1.252 \times 10^{-4}$ & $1.482 \times 10^{-5}$ 
 & $1.822 \times 10^{-6}$ & $2.268 \times 10^{-7}$
\\
\hline
err($\chi$) & $6.183 \times 10^{-2}$ & $8.600 \times 10^{-3}$ & $1.136 \times 10^{-3}$ 
 & $1.428 \times 10^{-4}$ & $1.814 \times 10^{-5}$
\\
\hline
\end{tabular}
\end{table}
\end{center}
\textbf{Remark.}\\
The numerical simulations show the robustness and the accuracy of the method both for fine 
and coarse meshes.
The differences in the errors between test A and B depend on the smoothness of the 
solutions; in the test B the maximum absolute value of the partial derivatives with respect to the coordinate $x$ is greater then in the first case. Moreover the different order of magnitude of the physical parameters introduces a stiffness in the set of partial differential equations.
\section{Spatial depending parameters}
\begin{figure}[ht]
	\centering
	\includegraphics[width=0.5\textwidth]{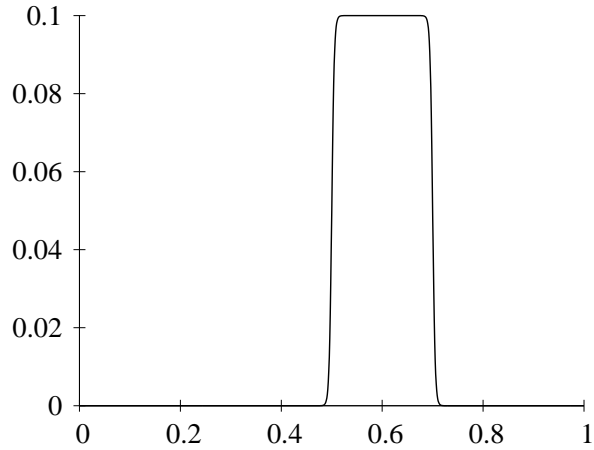}
	\caption{The function $\psi(x)$ in the interval $(0,1)$ for $h = 0.1$}
	\label{fig_psi}
\end{figure}

\begin{figure}[ht]
\centering
\includegraphics[width=0.49\textwidth]{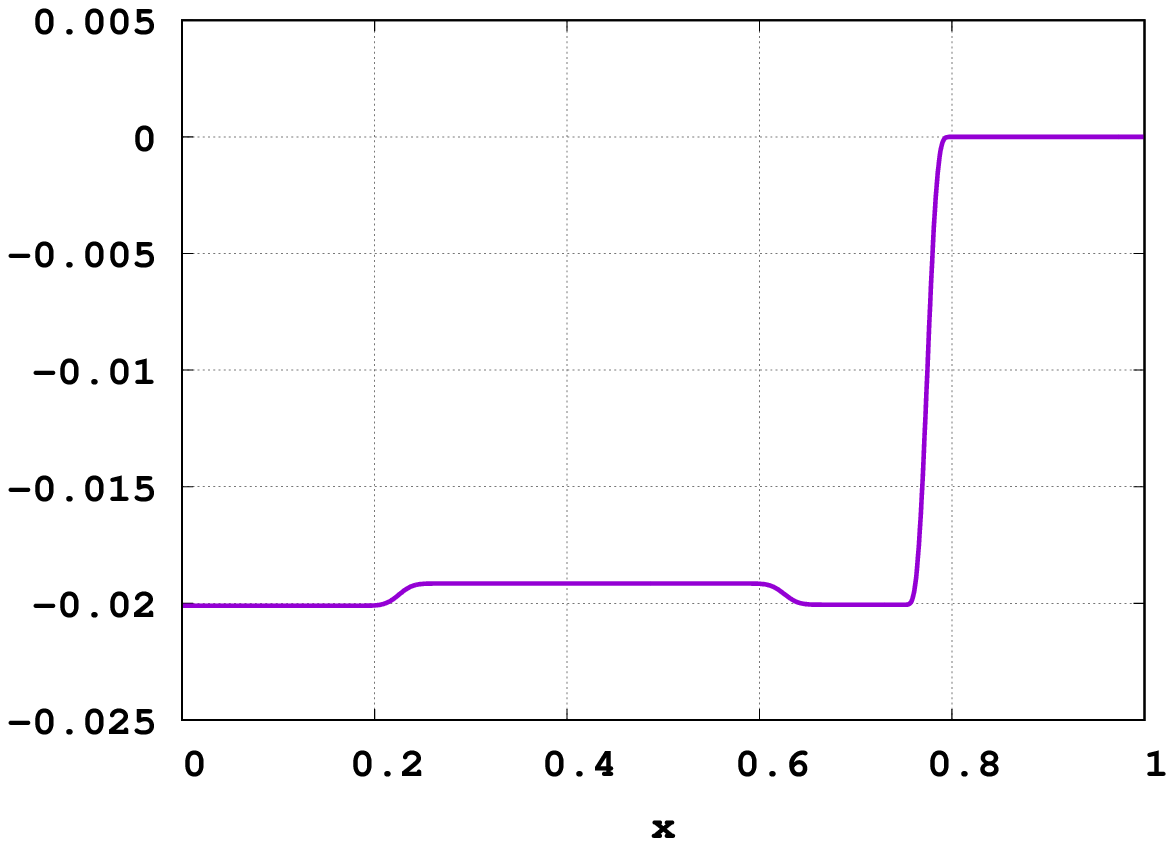}
\includegraphics[width=0.49\textwidth]{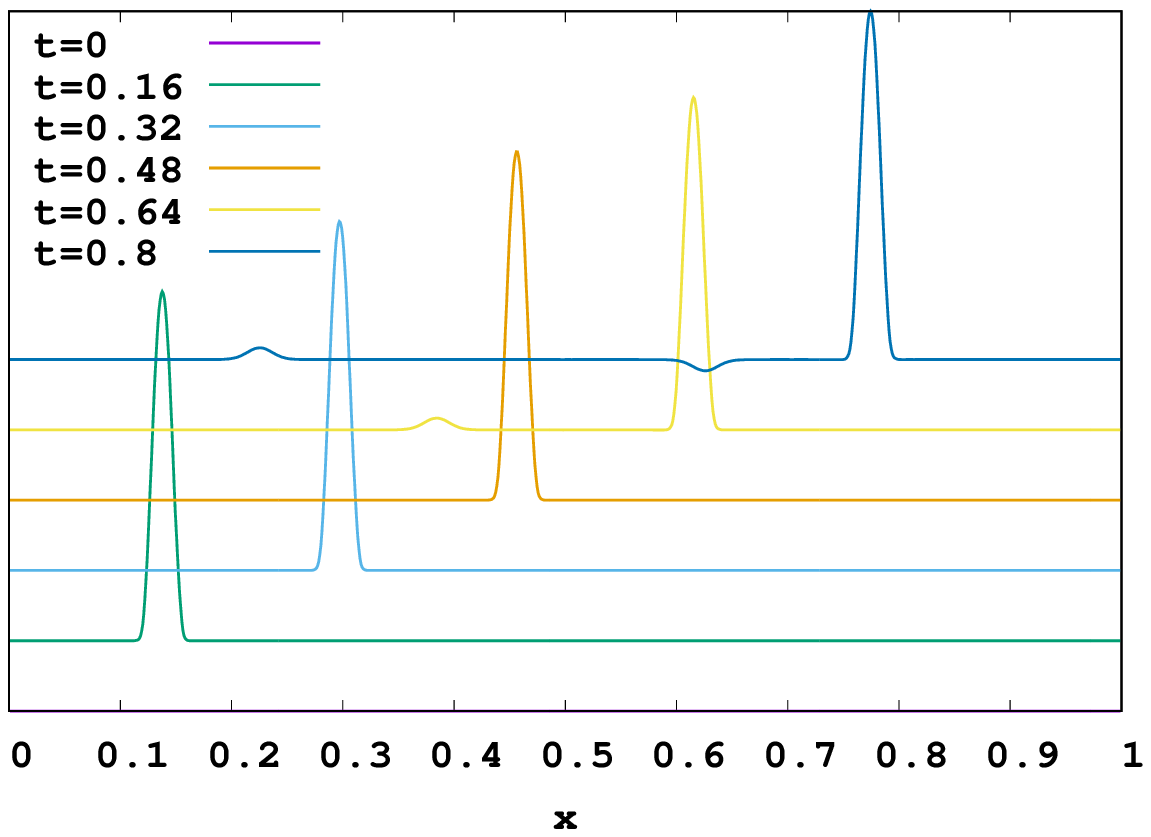}
\caption{The case $h = 0.1$. The function $u(t,x)$ at time $t = 0.8$ (left) and the function $u_{x}(t,x) +  \kappa \, t$ (right)}
\label{fig_ux01}
\end{figure}

\begin{figure}[ht]
\centering
\includegraphics[width=0.49\textwidth]{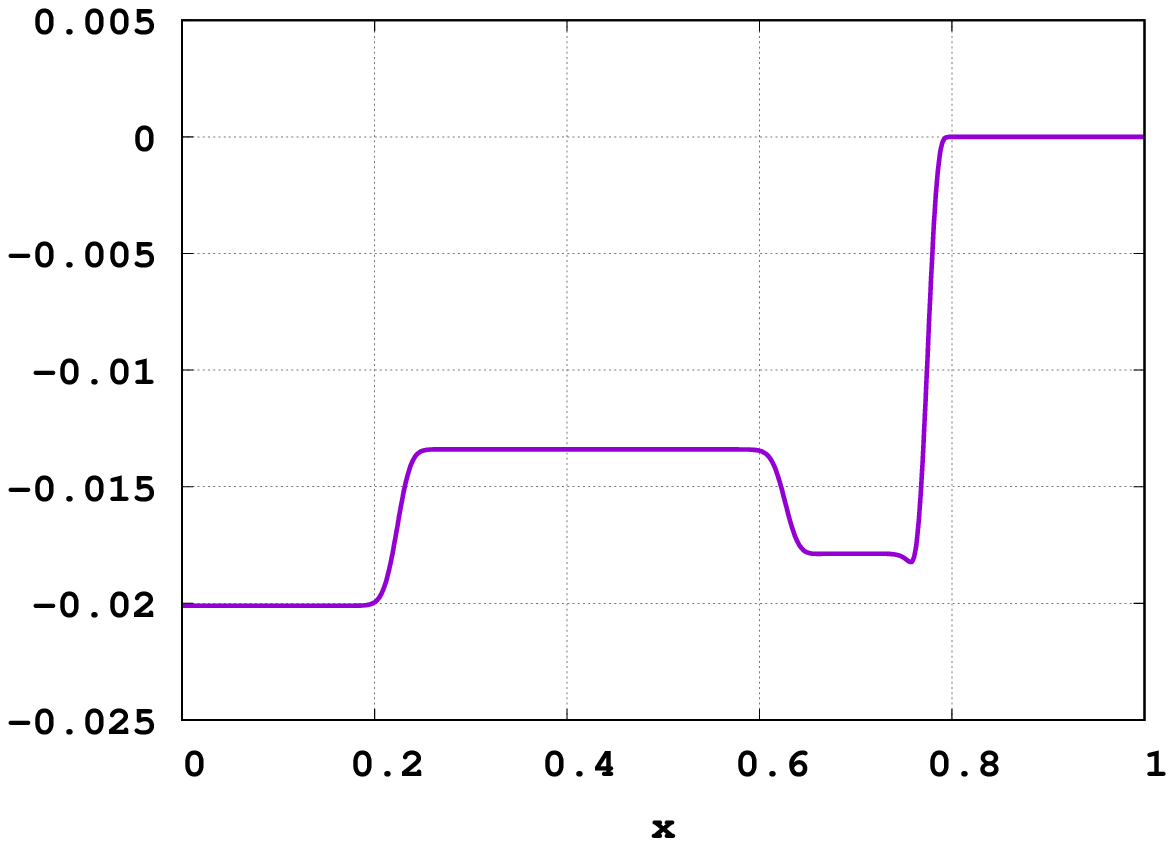}
\includegraphics[width=0.49\textwidth]{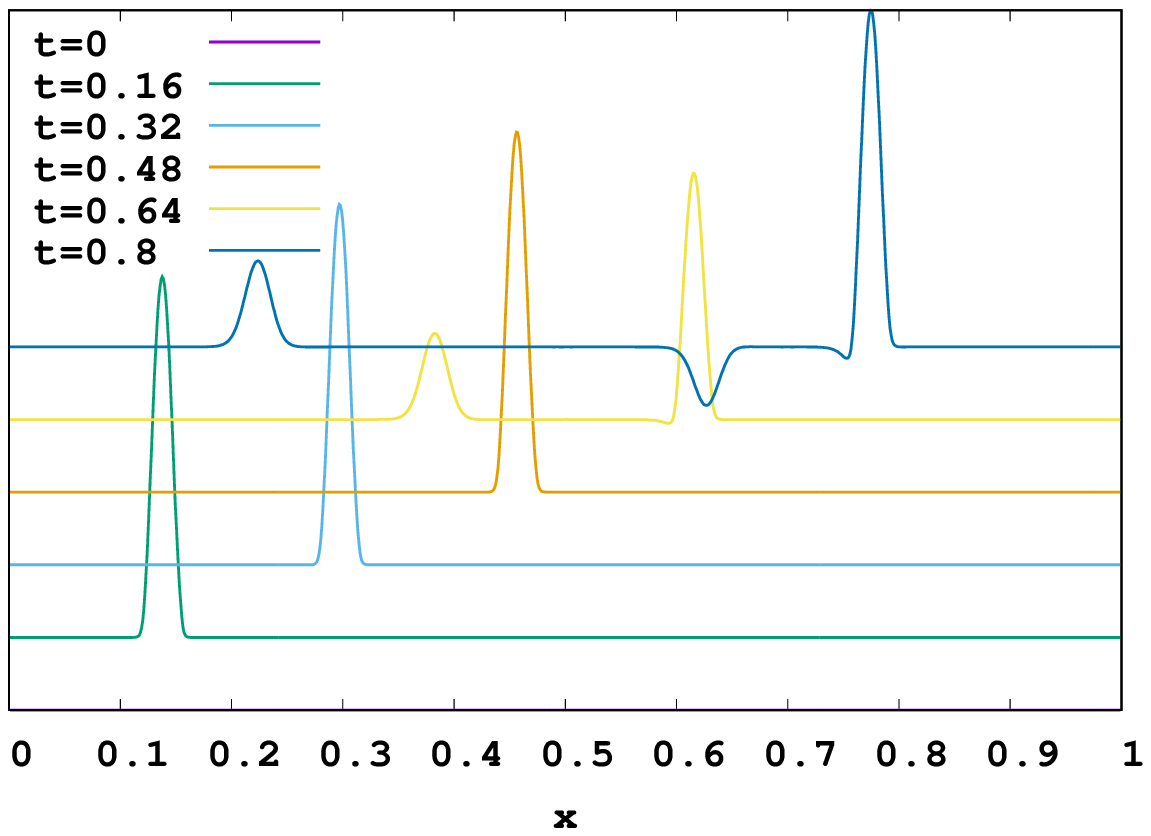}
\caption{The case $h = 1$. The function $u(t,x)$ at time $t = 0.8$ (left) and the function $u_{x}(t,x) +  \kappa \, t$ (right)}
\label{fig_ux1}
\end{figure}

When the physical parameters $\rho$, $I_{\mu}$, $\gamma$, $A$, $B$ and $C$ are not 
constant, but they are differentiable functions of the variable $x$, and we assume valid the definitions of the kinetic energy \eq{kinetic} and the potential energy \eq{potential}, then the Euler-Lagrange equations write
\begin{align}
\rho \, u_{tt} & =
\gamma \, u_{xx} + \dfrac{d \gamma}{d x} \, u_{x}  
+ A \, \chi_{x} + \dfrac{d A}{d x} \, \chi \sv
\\[7pt]
I_{\mu} \, \chi_{tt} & = 
C \, \chi_{xx} + \dfrac{d C}{d x} \, \chi_{x} - A \, u_{x} - B \, \chi \p
\end{align}
If we define
$$
\begin{array}{l}
a_{1}(x) = \dfrac{\gamma}{\rho} \sv \quad
b_{1}(x) = \dfrac{1}{\rho} \, \dfrac{d \gamma}{d x} \sv \quad
a_{2}(x) = \dfrac{A}{\rho} \sv \quad 
b_{2}(x) = \dfrac{1}{\rho} \, \dfrac{d A}{d x} \sv
\\[8pt]
a_{3}(x) = \dfrac{C}{I_{\mu}} \sv \quad
b_{3}(x) = \dfrac{1}{I_{\mu}} \, \dfrac{d C}{d x} \sv \quad 
a_{4}(x) = \dfrac{A}{I_{\mu}} \sv \quad 
a_{5}(x) = \dfrac{B}{I_{\mu}} \sv
\end{array}
$$
then the equations for the unknowns $u$ and $\chi$ becomes
\begin{align}
u_{tt} & =
a_{1}(x) \, u_{xx} + b_{1}(x) \, u_{x} + a_{2}(x) \, \chi_{x} + b_{2}(x) \, \chi \sv
\label{equ_tt_x}
\\[7pt]
\chi_{tt} & = 
a_{3}(x) \, \chi_{xx} + b_{3}(x) \, \chi_{x} - a_{4}(x) \, u_{x} - a_{5}(x) \, \chi \p
\label{eqX_tt_x}
\end{align}
Also in this case, we introduce new variables in order to make the system of partial 
differential equations suitable for a numerical integration.
It is possible to prove (see Appendix A) that Eqs.~\eq{equ_tt_x}-\eq{eqX_tt_x} are equivalent to the four partial differential equations 
\begin{align}
u_{t} & = \sqrt{a_{1}(x)} \, u_{x} - \varphi_{1}(x) \, u - \varphi_{2}(x) \, \chi + v \sv
\label{eq_ux}
\\
\chi_{t} & = \sqrt{a_{3}(x)} \, \chi_{x} - \varphi_{3}(x) \, u - \varphi_{4}(x) \, \chi + w 
\sv
\\
v_{t} & = - \sqrt{a_{1}(x)} \, v_{x}
+ \phi_{1}(x) \, u + \phi_{2}(x) \, \chi + \varphi_{1}(x) \, v + \varphi_{2}(x) \, w \sv
\\
w_{t} & = - \sqrt{a_{3}(x)} \, w_{x} + \phi_{3}(x) \, u + \phi_{4}(x) \, \chi
+ \varphi_{3}(x) \, v + \varphi_{4}(x) \, w \sv
\end{align}
where
\begin{align*}
&
\phi_{1}(x) = \sqrt{a_{1}(x)} \: \dfrac{d \varphi_{1}}{d x} -
\left[ \varphi_{1}(x) \right]^{2} -  \varphi_{2}(x) \, \varphi_{3}(x) \sv
\\[7pt]
&
\phi_{2}(x) = \sqrt{a_{1}(x)} \: \dfrac{d \varphi_{2}}{d x}  + b_{2}(x)
- \varphi_{1}(x) \, \varphi_{2}(x) - \varphi_{2}(x) \, \varphi_{4}(x) \sv
\\[7pt]
&
\phi_{3}(x) = \sqrt{a_{3}(x)} \: \dfrac{d \varphi_{3}}{d x} -
\varphi_{3}(x) \, \varphi_{1}(x) - \varphi_{4}(x) \, \varphi_{3}(x) \sv
\\[7pt]
&
\phi_{4}(x) = \sqrt{a_{3}(x)} \: \dfrac{d \varphi_{4}}{d x} - a_{5}(x) - 
\varphi_{3}(x) \, \varphi_{2}(x) - \left[ \varphi_{4}(x) \right]^{2} ,
\\[7pt]
&
\varphi_{1}(x) = \dfrac{-1}{2 \, \rho(x)} 
\dfrac{d \mbox{ }}{d x} \left[ \rho(x) \, \sqrt{a_{1}(x)} \right] 
, \quad
\varphi_{2}(x) = \dfrac{\mbox{} - a_{2}(x)}{\sqrt{a_{1}(x)} + \sqrt{a_{3}(x)}} \sv
\\[7pt]
&
\varphi_{3}(x) = \dfrac{a_{4}(x)}{\sqrt{a_{1}(x)} + \sqrt{a_{3}(x)}}
\sv \quad
\varphi_{4}(x) = \dfrac{-1}{2 \, I_{\mu}(x)} \dfrac{d \mbox{ }}{d x} 
\left[ \sqrt{a_{3}(x)} \, I_{\mu}(x) \right] .
\end{align*}
We show a simple numerical example, where the numerical scheme is the same of the cases 
described in the previous section. The spatial domain is the interval $[0,1]$.
The assume that
\begin{align*}
&
\gamma(x) = \gamma^{*} \, (1 + \psi(x)) \sv \quad
A(x) = A^{*} \, (1 + \psi(x)) \sv \quad
B(x) = B^{*} \, (1 + \psi(x)) 
\\
&
\rho(x) = \rho^{*} \, (1 + \psi(x)) \sv \quad
C(x) = C^{*} \, (1 + \psi(x)) \sv \quad
I_{\mu}(x) = I_{\mu}^{*} \, (1 + \psi(x)) \sv
\end{align*}
with
$$
\rho^{*} = 1 \sv \quad I_\mu^{*} = 1 \sv \quad
\gamma^{*} = 0.99 \sv \quad A^{*} = - \, 0.01 \sv \quad B^{*} = 10 \sv \quad C^{*} = 1 
\sv
$$
and (see, Figure~\ref{fig_psi}) 
$$
\psi(x) = h \left[ \dfrac{1}{1 + \exp(400 \, (0.5 - x))} -
\dfrac{1}{1 + \exp(400 \, (0.7 - x))} \right] ,
$$
where $h$ is a parameter. In our simulations we have chosen $h = 0.1$ and $h = 1$.
We point out that the function $\psi$ is not smooth. The model simulates an approximation 
of a continuum, which consists of two different materials.
Following Ref.~\cite{Berezovski2016b} we assume that the continuum is at rest at the initial time, that is
$$
u(0,x) = 0 \sv \quad \chi(0,x) = 0 \sv \quad 
v(0,x) = 0 \sv \quad w(0,x) = 0 \p
$$
The boundary conditions must be simulated an excitation of the strain at
$x = 0$ for an short time period; so we must assign $u_{x}(t,0^{-})$.
This is possible, by observing that only Eq.~\eq{eq_ux} contains the term $u_{x}$.
Now, if $u_{x}(t,0^{-}) = \varepsilon(t)$, then Eq.~\eq{eq_ux} gives
\begin{equation}
v(t,0^{-}) = u_{t}(t,0) - \sqrt{a_{1}(0)} \, \varepsilon(t) + \varphi_{1}(0) \, u(t,0) + \varphi_{2}(0) \, \chi(t,0) \sv
\label{bc_v}
\end{equation}
which is the boundary condition for the unknown $v$ at time $t$.
The right hand side of \eq{bc_v} can be evaluated, since it is known at the previous time step. The partial derivative of $u$ with respect to the time is approximated by means of a simple finite difference.
The other boundary conditions \eq{bc} are assumed null at every time. 
So, when we find the approximated solution at each time step by means of numerical integration.
We have chosen
$$
\varepsilon(t) =
\left\{
\begin{array}{ll}
\dfrac{1 + \cos(\pi \, (1 - 50 \, t))}{2} & \mbox{ if } 0 \leq t \leq 0.04
\\
0 & \mbox{ otherwise}
\end{array}
\right.
$$
in the simulation.
 \\
We have found the solutions (see Figure~\ref{fig_ux01} and \ref{fig_ux1}) for $h = 0.1$ and $h = 1$, by using $1024$ grid points in the spatial interval $[0,1]$, and we verify a nice accuracy by changing the number of grid points. In the figure the partial derivative of $u(t,x)$ with respect to $x$ is shown at different times.
We plot $u_{x}(t,x) + \kappa \, t$, where $\kappa$ is a positive constant, instead of $u_{x}(t,x)$. The fictitious shift was introduced only to make clear the figure.
The reflected waves, which arise from the inhomogeneity, are evident.
\section{Conclusions}
In this paper a class of exact  solutions of a one dimensional Mindlin model to describe 
linear elastic behaviour is obtained.
The main assumption of this paper concerns the potential energy, which is assumed  
strictly positive definite.
Due to this property the solutions of the model have a wave behaviour.
We find exact solutions to test the accuracy of the proposed numerical method. It is based 
on a weighted essentially non-oscillatory finite difference scheme, coupled by a total 
variation diminishing Runge-Kutta method.
The results obtained with the numerical methods were matched to the exact analytical 
solutions, and they agree very well. 
In this way we also verify the robustness and the accuracy of numerical scheme both for fine and coarse meshes.
Also in the case when the physical parameters are not constant and smooth, the numerical scheme seems to give accurate solutions to the equations.
\section*{Acknowledgments}
The first author was partially supported by the italian FIR project "Innovative techniques 
in computational mechanics based on high continuity interpolation for the integrated design of advanced structures" (Principal Investigator: Massimo Cuomo).
\section*{Appendix A}
Let us consider Eqs.~\eq{equ_tt_x}-\eq{eqX_tt_x}.
Firstly we define the new unknowns
$$
\alpha = u_{t} - \sqrt{a_{1}(x)} \, u_{x}  \quad \mbox{and} \quad
\beta = \chi_{t} - \sqrt{a_{3}(x)} \, \chi_{x} \p
$$
So Eqs.~\eq{equ_tt_x}-\eq{eqX_tt_x} are equivalent to the system
\begin{align}
u_{t} & = \sqrt{a_{1}(x)} \, u_{x} + \alpha \sv
\label{equx}
\\
\chi_{t} & = \sqrt{a_{3}(x)} \, \chi_{x} + \beta \sv
\\
\alpha_{t} & = \mbox{} - \sqrt{a_{1}(x)} \, \alpha_{x} + \left[ b_{1}(x) -
\dfrac{1}{2} \, \dfrac{d a_{1}(x)}{d x} \right] u_{x}
 + a_{2}(x) \, \chi_{x} + b_{2}(x) \, \chi \sv
\\
\beta_{t} & = \mbox{} - \sqrt{a_{3}(x)} \, \beta_{x} + \left[ b_{3}(x) -
\dfrac{1}{2} \, \dfrac{d a_{3}(x)}{d x} \right] \chi_{x} 
 - a_{4}(x) \, u_{x} - a_{5}(x) \, \chi \p
 \label{eqbx}
\end{align}
Since
\begin{align*}
&
b_{1}(x) - \dfrac{1}{2} \, \dfrac{d a_{1}(x)}{d x} = 
\dfrac{1}{\rho(x)} \, \dfrac{d \mbox{ }}{d x} \left[ \rho(x) \, a_{1}(x) 
\right] - \dfrac{1}{2} \, \dfrac{d a_{1}(x)}{d x} =
\dfrac{a_{1}(x)}{\rho(x)} \, \dfrac{d \rho(x)}{d x} + 
\dfrac{1}{2} \, \dfrac{d a_{1}(x)}{d x} \sv
\\
&
b_{3}(x) - \dfrac{1}{2} \, \dfrac{d a_{3}(x)}{d x} =
\dfrac{1}{I_{\mu}(x)} \, \dfrac{d \mbox{ }}{d x} \left[ I_{\mu}(x) \, a_{3}(x) 
\right] - \dfrac{1}{2} \, \dfrac{d a_{3}(x)}{d x} =
\dfrac{a_{3}(x)}{I_{\mu}(x)} \, \dfrac{d I_{\mu}(x)}{d x} +
\dfrac{1}{2} \, \dfrac{d a_{3}(x)}{d x} \sv
\end{align*}
then Eqs.~\eq{equx}-\eq{eqbx} write
\begin{align*}
u_{t} & = \sqrt{a_{1}(x)} \, u_{x} + \alpha \sv
\\
\chi_{t} & = \sqrt{a_{3}(x)} \, \chi_{x} + \beta \sv
\\
\alpha_{t} & = \mbox{} - \sqrt{a_{1}(x)} \, \alpha_{x} + \left[ \dfrac{a_{1}(x)}{\rho(x)} \, \dfrac{d \rho(x)}{d x} + \dfrac{1}{2} \, \dfrac{d a_{1}(x)}{d x} 
\right] u_{x} + a_{2}(x) \, \chi_{x} + b_{2}(x) \, \chi \sv
\\
\beta_{t} & = \mbox{} - \sqrt{a_{3}(x)} \, \beta_{x} + \left[ \dfrac{a_{3}(x)}{I_{\mu}(x)} 
\, \dfrac{d I_{\mu}(x)}{d x} + \dfrac{1}{2} \, \dfrac{d a_{3}(x)}{d x} \right] \chi_{x} 
- a_{4}(x) \, u_{x} - a_{5}(x) \, \chi \p
\end{align*} 
Now let's introduce the new unknowns $v$ and $w$ to replace $\alpha$ and $\beta$, by means 
of the relationship
$$
v = \alpha + \varphi_{1}(x) \, u + \varphi_{2}(x) \, \chi \sv 
\quad
w = \beta + \varphi_{3}(x) \, u + \varphi_{4}(x) \, \chi \sv
$$
where the functions $\varphi_{i}$ must be chosen opportunely.
It is a simple matter to verify that the new equations are
\begin{align*}
u_{t} & = \sqrt{a_{1}(x)} \, u_{x} - \varphi_{1}(x) \, u - \varphi_{2}(x) \, \chi + v \sv 
\\
\chi_{t} & = \sqrt{a_{3}(x)} \, \chi_{x} - \varphi_{3}(x) \, u - \varphi_{4}(x) \, \chi + w \sv
\\
v_{t} & = \left[ \sqrt{a_{1}(x)} \, \varphi_{1}(x) + \dfrac{a_{1}(x)}{\rho(x)} \, 
\dfrac{d \rho(x)}{d x} + \dfrac{1}{2} \, \dfrac{d a_{1}(x)}{d x} 
+ \varphi_{1}(x) \, \sqrt{a_{1}(x)} \right] u_{x} 
\\
& \quad
+ \left[ \sqrt{a_{1}(x)} \, \varphi_{2}(x) + a_{2}(x) + \varphi_{2}(x) \, \sqrt{a_{3}(x)}
\right] \chi_{x} - \sqrt{a_{1}(x)} \, v_{x} 
\\
& \quad 
+ \left[ \sqrt{a_{1}(x)} \, \dfrac{d \varphi_{1}}{d x} -
\left[ \varphi_{1}(x) \right]^{2} -  \varphi_{2}(x) \, \varphi_{3}(x)
\right] u
\\
& \quad 
+ \left[ \sqrt{a_{1}(x)} \, \dfrac{d \varphi_{2}}{d x}  + b_{2}(x)
- \varphi_{1}(x) \, \varphi_{2}(x) - \varphi_{2}(x) \, \varphi_{4}(x) \right] \chi 
+ \varphi_{1}(x) \, v + \varphi_{2}(x) \, w \sv
\\
w_{t} & =
\left[ \sqrt{a_{3}(x)} \, \varphi_{3}(x) - a_{4}(x) + \varphi_{3}(x)  \sqrt{a_{1}(x)}
\right] u_{x}
\\
& \quad + \left[ \sqrt{a_{3}(x)} \, \varphi_{4}(x) +
\dfrac{a_{3}(x)}{I_{\mu}(x)} \, \dfrac{d I_{\mu}(x)}{d x} +
\dfrac{1}{2} \, \dfrac{d a_{3}(x)}{d x}
+ \varphi_{4}(x) \sqrt{a_{3}(x)}  \right] \chi_{x} - \sqrt{a_{3}(x)} \, w_{x}
\\
& \quad +
\left[ \sqrt{a_{3}(x)} \, \dfrac{d \varphi_{3}}{d x} -
\varphi_{3}(x) \, \varphi_{1}(x) - \varphi_{4}(x) \, \varphi_{3}(x)  \right] u 
\\
& \quad +
\left[  \sqrt{a_{3}(x)} \, \dfrac{d \varphi_{4}}{d x} - a_{5}(x) - 
\varphi_{3}(x) \, \varphi_{2}(x) - \left[ \varphi_{4}(x) \right]^{2} \right] \chi
+ \varphi_{3}(x) \, v + \varphi_{4}(x) \, w \p
\end{align*}
We choose the functions $\varphi_{i}$ to make null the coefficients of $u_{x}$ and 
$\chi_{x}$ in the third and fourth equation.
We derive immediately
\begin{align*}
\varphi_{1}(x) & = \dfrac{-1}{2 \, \sqrt{a_{1}(x)}} \left[
\dfrac{a_{1}(x)}{\rho(x)} \, \dfrac{d \rho(x)}{d x} + 
\dfrac{1}{2} \, \dfrac{d a_{1}(x)}{d x} \right]
= \dfrac{-1}{2 \, \rho(x)} 
\dfrac{d \mbox{ }}{d x} \left[ \rho(x) \, \sqrt{a_{1}(x)} \right] ,
\\
\varphi_{2}(x) & = \dfrac{\mbox{} - a_{2}(x)}{\sqrt{a_{1}(x)} + \sqrt{a_{3}(x)}} \sv
\\
\varphi_{3}(x) & = \dfrac{a_{4}(x)}{\sqrt{a_{1}(x)} + \sqrt{a_{3}(x)}} \sv
\\
\varphi_{4}(x) & = \dfrac{-1}{2 \, \sqrt{a_{3}(x)}} \left[ \dfrac{a_{3}(x)}{I_{\mu}(x)} \, 
\dfrac{d I_{\mu}(x)}{d x} +
\dfrac{1}{2} \, \dfrac{d a_{3}(x)}{d x} \right] 
=
\dfrac{-1}{2 \, I_{\mu}(x)} \dfrac{d \mbox{ }}{d x} 
\left[ \sqrt{a_{3}(x)} \, I_{\mu}(x) \right] .
\end{align*}
Therefore the set of equations reduces to
\begin{align*}
u_{t} & = \sqrt{a_{1}(x)} \, u_{x} - \varphi_{1}(x) \, u - \varphi_{2}(x) \, \chi + v \sv
\\
\chi_{t} & = \sqrt{a_{3}(x)} \, \chi_{x} - \varphi_{3}(x) \, u - \varphi_{4}(x) \, \chi + w 
\sv 
\\
v_{t} & = - \sqrt{a_{1}(x)} \, v_{x}
+ \left[ \sqrt{a_{1}(x)} \, \dfrac{d \varphi_{1}}{d x} -
\left[ \varphi_{1}(x) \right]^{2} -  \varphi_{2}(x) \, \varphi_{3}(x)
\right] u 
\\
&
\quad
+ \left[ \sqrt{a_{1}(x)} \, \dfrac{d \varphi_{2}}{d x}  + b_{2}(x)
- \varphi_{1}(x) \, \varphi_{2}(x) - \varphi_{2}(x) \, \varphi_{4}(x) \right] \chi
+ \varphi_{1}(x) \, v + \varphi_{2}(x) \, w \sv
\\
w_{t} & = - \sqrt{a_{3}(x)} \, w_{x} +
\left[ \sqrt{a_{3}(x)} \, \dfrac{d \varphi_{3}}{d x} -
\varphi_{3}(x) \, \varphi_{1}(x) - \varphi_{4}(x) \, \varphi_{3}(x)  \right] u 
\\
&
\quad +
\left[  \sqrt{a_{3}(x)} \, \dfrac{d \varphi_{4}}{d x} - a_{5}(x) - 
\varphi_{3}(x) \, \varphi_{2}(x) - \left[ \varphi_{4}(x) \right]^{2} \right] \chi
+ \varphi_{3}(x) \, v + \varphi_{4}(x) \, w  \p
\end{align*}
%
%
%
%
%

%
\end{document}